\documentclass[12pt]{article}
\usepackage[utf8]{inputenc}
\usepackage[english%, russian
]{babel}
\usepackage {amsmath, amssymb, color}
\def\ep{\varepsilon}

\newtheorem{Theorem}{Theorem}
\newtheorem{Lemma}{Lemma}
\newtheorem{Remark}[Theorem]{Remark}

\begin{document}

\title %{\vspace{-10mm}
{Some lemmata on the perturbation\\ of the spectrum}

\author{Alexander I. Nazarov\footnote{
St.Petersburg Department of Steklov Institute, Fontanka 27, St.Petersburg, 191023, Russia, 
and St.Petersburg State University, 
Universitetskii pr. 28, St.Petersburg, 198504, Russia. E-mail: al.il.nazarov@gmail.com.
}
}
\date{}

\maketitle

\begin{abstract}
\footnotesize
We give some sufficient conditions for preserving of the second term in the spectral asymptotics of a compact operator under the perturbation of the metrics in the Hilbert space.
\medskip
\end{abstract}

It is well known, see, e.g., \cite[Lemma 1.16]{BS74}, that the one-term power-type spectral asymptotics of a compact operator in the Hilbert space does not change under compact perturbation 
of the metrics of the space. The problem of preserving of the two-term asymptotics is much more sensitive and complicated. Here we give some sufficient conditions for this. These results 
can be applied in the spectral analysis of some integro-differential operators arising in the theory of Gaussian random processes, see \cite{Naz19}.
\medskip

In what follows we denote by $c$ any absolute constant.

\begin{Lemma}\label{norma}
 Let ${\cal K}$ and ${\cal B}$ be self-adjoint compact operators in the Hilbert space ${\cal H}$. Suppose that ${\cal K}$ and ${\cal I}+{\cal B}$ are positive. Denote by $\lambda_n$ 
 the eigenvalues of ${\cal K}$ enumerated in the decreasing order taking into account the multiplicities, and by $h_n$ corresponding normalized eigenfunctions. Finally, suppose that
\begin{equation}\label{asymp}
\lambda_n=\big(an+b+O(n^{-\delta})\big)^{-B}, \qquad \|{\cal B}h_n\|_{\cal H}\le cn^{-(1+\delta)},
\end{equation}
 as $n\to\infty$, where $a, B,\delta>0$, $b\in\mathbb{R}$. 
 Then the eigenvalues $\mbox{\boldmath$\lambda$}_n$ of generalized eigenproblem
\begin{equation}\label{gen-eigen}
{\cal K}{\bf h}_n=\mbox{\boldmath$\lambda$}_n\big({\bf h}_n+{\cal B}{\bf h}_n\big)
\end{equation}
 have the same two-term asymptotics as $n\to\infty$:
 \begin{equation}\label{asymp1}
\mbox{\boldmath$\lambda$}_n=\big(an+b+O(n^{-\delta})\big)^{-B}. 
\end{equation}
\end{Lemma}

{\bf Proof.} We introduce new scalar products in ${\cal H}$:
 \begin{equation*}
\langle h,g\rangle:=(h+{\cal B}h,g)_{\cal H}.
 \end{equation*}
It is easy to see that corresponding norm $|\!|\!|h|\!|\!|:=\langle h,h\rangle^{\frac 12}$ is equivalent to original one. Denote by $\mathbb{H}$ the space ${\cal H}$ with new scalar product. 
Then the sesquilinear form $({\cal K}h,g)_{\cal H}$ generates a compact positive self-adjoint operator $\mathbb{B}$ such that
 \begin{equation*}
\langle \mathbb{B}h,g\rangle=({\cal K}h,g)_{\cal H},\qquad h,g\in\mathbb{H},
 \end{equation*}
 and the generalized eigenproblem (\ref{gen-eigen}) is reduced to the standard eigenproblem for the operator $\mathbb{B}$ in $\mathbb{H}$. 

Recall some elementary facts from the theory of spectral measure, see Ch. 5 in \cite{BS10}. 
The spectral measure $d{\cal E}(t)$ associated with $\mathbb{B}$ generates the family of scalar measures
$$
de_h(t):=\langle d{\cal E}(t)h,h\rangle, \qquad h\in \mathbb{H}.
$$
Moreover, the following obvious formulae hold for arbitrary $h\in \mathbb{H}$:
\begin{equation*}%\label{spectral_int}
|\!|\!|h|\!|\!|^2=\int\limits_{\mathbb R} de_h(t), \qquad
|\!|\!| \mathbb{B} h-\lambda h|\!|\!|^2= \int\limits_{\mathbb R} (t-\lambda)^2de_h(t).
\end{equation*}

If we assume that an interval $\Delta=(\lambda-\delta,\lambda+\delta)$ is free of the spectrum of $\mathbb{B}$ then we have for any $h\in \mathbb{H}$
\begin{equation}\label{free_segment}
|\!|\!|\mathbb{B} h-\lambda h|\!|\!|^2= \int\limits_{\mathbb R\setminus\Delta} (t-\lambda)^2de_h(t) \ge \delta^2\int\limits_{\mathbb R\setminus\Delta} de_h(t)=\delta^2\,|\!|\!|h|\!|\!|^2.
\end{equation}

Now we set $\lambda=\lambda_n$, $h=h_n$. For any $g\in \mathbb{H}$ we have
\begin{equation*}
\aligned
|\langle\mathbb{B} h_n-\lambda_n h_n,g\rangle| & =|({\cal K}h_n,g)_{\cal H}-\lambda_n(h_n+{\cal B}h_n,g)_{\cal H}|\\
& =\lambda_n|({\cal B}h_n,g)_{\cal H}|
\le \lambda_n \|{\cal B}h_n\|_{\cal H}\|g\|_{\cal H}\le\frac {c_1}{n^{1+\delta}}\,\lambda_n|\!|\!|g|\!|\!|,
\endaligned
\end{equation*}
and therefore 
\begin{equation*}
|\!|\!|\mathbb{B} h_n-\lambda_n h_n|\!|\!|\le\frac {c_1}{n^{1+\delta}}\,\lambda_n\le\frac {c_2}{n^{1+\delta}}\,\lambda_n|\!|\!|h_n|\!|\!|.
\end{equation*}
Comparing this inequality with (\ref{free_segment}) we see that the interval 
\begin{equation*}
\Delta_n=\big(\lambda_n(1-c_2n^{-(1+\delta)}),\lambda_n(1+c_2n^{-(1+\delta)})\big)
\end{equation*}
contains an eigenvalue $\mbox{\boldmath$\lambda$}$ of the generalized eigenproblem (\ref{gen-eigen}). 

By (\ref{asymp}), intervals $\Delta_n$ and $\Delta_{n+1}$ do not intersect for $n$ sufficiently large. Repeating previous argument for $\ep{\cal B}$ instead of ${\cal B}$, $\ep\in[0,1]$, 
we notice that the eigenvalues depend continuously on $\varepsilon$ and conclude that for large $n$ the interval $\Delta_n$ contains just $\mbox{\boldmath$\lambda$}_n$. This yields 
(\ref{asymp1}).
\hfill$\square$\medskip

This result is quite simple but the assumption (\ref{asymp}) is very restrictive. The following theorem gives a ``more pointwise'' condition which is, however, globally weaker.

\begin{Theorem}
 In Lemma \ref{norma}, suppose that instead of (\ref{asymp}) the following relations hold:
\begin{equation*}%\label{asymp2}
\lambda_n=\big(an+b+O(n^{-\delta})\big)^{-B}, \qquad |({\cal B}h_n,h_m)_{\cal H}|\le c(mn)^{-\frac {1+\delta}2}.
\end{equation*}
 Then (\ref{asymp1}) also holds.
 \end{Theorem}

{\bf Proof.} First, we notice that we can write ${\cal B}={\cal B}_++{\cal B}_-$, where ${\cal B}_+\ge0$ and ${\cal B}_-\le0$. By the min-max principle (see, e.g., \cite[Appendix 1]{BS74}, we have 
$\mbox{\boldmath$\lambda$}_n^+\le\mbox{\boldmath$\lambda$}_n\le\mbox{\boldmath$\lambda$}_n^-$, where $\mbox{\boldmath$\lambda$}_n^{\pm}$ are eigenvalues of the problems
\begin{equation*}
{\cal K}{\bf h}_n^+=\mbox{\boldmath$\lambda$}_n^+\big({\bf h}_n^++{\cal B}_+{\bf h}_n^+\big); \qquad {\cal K}{\bf h}_n^-=\mbox{\boldmath$\lambda$}_n^-\big({\bf h}_n^-+{\cal B}_-{\bf h}_n^-\big).
\end{equation*}
So, it suffices to consider two cases: positive ${\cal B}$ and negative ${\cal B}$.
\medskip

{\bf 1}. Let ${\cal B}$ be positive. Then evidently $\mbox{\boldmath$\lambda$}_n\le \lambda_n$. On the other hand, the min-max principle gives $\mbox{\boldmath$\lambda$}_{n}\ge \widehat\lambda_n$, 
where $\widehat\lambda_k$ are the eigenvalues of generalized (finite-dimensional) eigenproblem
\begin{equation*}%\label{hat-eigen-1}
\widehat P_n{\cal K}\widehat P_n\widehat h_k=\widehat\lambda_k\big(\widehat h_k+\widehat P_n{\cal B}\widehat P_n\widehat h_k\big),
\end{equation*}
and $\widehat P_n$ is the orthoprojector onto the span of $\widehat{\cal H}_n=\text{Span}\{h_k\}$, $k\le n$.

Let $\widehat x\in \widehat{\cal H}_n$ be the minimizer of the Rayleigh quotient  
\begin{equation}\label{Rayleigh}
J(x)=\frac {({\cal K}x,x)_{\cal H}}{(x,x)_{\cal H}+({\cal B}x,x)_{\cal H}}
\end{equation}
over $\widehat{\cal H}_n$. We derive for $k< n$
\begin{equation*}%\label{deriv}
0=\frac 12\,J'(\widehat x;h_k)=\frac {({\cal K}\widehat x,h_k)_{\cal H}-J(\widehat x)\cdot\big((\widehat x,h_k)_{\cal H}+({\cal B}\widehat x,h_k)_{\cal H})\big)}{(\widehat x,\widehat x)_{\cal H}
+({\cal B}\widehat x,\widehat x)_{\cal H}}.
\end{equation*}
Therefore,
\begin{equation*}
0=(\lambda_k-J(\widehat x))\cdot(\widehat x,h_k)_{\cal H}-J(\widehat x)({\cal B}\widehat x,h_k)_{\cal H},
\end{equation*}
i.e.
\begin{equation}\label{J'=0}
(\widehat x,h_k)_{\cal H}=\frac {J(\widehat x)}{\lambda_k-J(\widehat x)}\,({\cal B}\widehat x,h_k)_{\cal H}.
\end{equation}

Since $J(h_n)\le \lambda_n$, we have $J(\widehat x)\le\lambda_n$. So, for any $k<n$
\begin{equation*}%\label{ak}
\widehat a_k:=|(\widehat x,h_k)_{\cal H}|\le\frac {\lambda_n}{\lambda_k-\lambda_n}\,\sum\limits_{m=1}^n \widehat a_m\cdot \frac c{(km)^{\frac {1+\delta}2}}.
\end{equation*}
%Without loss of generality we set $\widehat a_n=1$. 
This implies 
%\begin{equation*}
%\widehat A:=\sum\limits_{k=1}^n \frac {\widehat a_k}{k^{\frac {1+\delta}2}}= \sum\limits_{k=1}^{n-1} \frac {\widehat a_k}{k^{\frac {1+\delta}2}}+\frac {\widehat a_n}{n^{\frac {1+\delta}2}}. 
%\end{equation*}
%Therefore, (\ref{ak}) gives
\begin{equation}\label{A}
\widehat A:=\sum\limits_{k=1}^n \frac {\widehat a_k}{k^{\frac {1+\delta}2}}\le \widehat A\,\sum\limits_{k=1}^{n-1}\frac {\lambda_n}{\lambda_k-\lambda_n}\cdot \frac c{k^{1+\delta}}
+\frac {\widehat a_n}{n^{\frac {1+\delta}2}}=:\widehat A\widehat {\mathfrak C}+\frac {\widehat a_n}{n^{\frac {1+\delta}2}}.
\end{equation}

Notice that for $k<n$
\begin{equation*}
\frac {\lambda_n}{\lambda_k-\lambda_n}\le c\,\frac {(k/n)^B}{1-(k/n)^B}+\frac c{n^{\min\{1,B\}}}.
\end{equation*}
So, 
\begin{equation*}
\aligned
\widehat {\mathfrak C}\le &\ \frac c{n^{1+\delta}}\,\sum\limits_{k=1}^{n-1} \frac {(k/n)^{B-1-\delta}}{1-(k/n)^B}+\frac c{n^{\min\{1,B\}}}\\
\le &\ \frac c{n^{\delta}}\int\limits_{\frac 1n}^{1-\frac 1n} \frac {t^{B-1-\delta}}{1-t^B}\,dt+\frac c{n^{\min\{1,B\}}}\le \frac {c\log(n)}{n^{\min\{1,\delta,B\}}},
\endaligned
\end{equation*}
and for $n$ sufficiently large (\ref{A}) gives $\widehat A\le c\widehat a_nn^{-\frac {1+\delta}2}$.

Now we calculate
\begin{equation*}
J(\widehat x)\ge\frac {\sum\limits_{k=1}^n \lambda_k \widehat a_k^2}{\sum\limits_{k=1}^n \widehat a_k^2+c\sum\limits_{k,m=1}^n\dfrac {\widehat a_k\widehat a_m}{(km)^{\frac {1+\delta}2}}}
\ge \lambda_n(1-cn^{-(1+\delta)}),
\end{equation*}
and the statement follows.
\medskip

{\bf 2}. Let ${\cal B}$ be negative. Then evidently $\mbox{\boldmath$\lambda$}_n\ge \lambda_n$. On the other hand, the min-max principle gives $\mbox{\boldmath$\lambda$}_{n+k-1}\le \widetilde\lambda_k$, 
where $\widetilde\lambda_k$ are the eigenvalues of generalized eigenproblem
\begin{equation*}%\label{hat-eigen}
\widetilde P_n{\cal K}\widetilde P_n\widetilde h_k=\widetilde\lambda_k\big(\widetilde h_k+\widetilde P_n{\cal B}\widetilde P_n\widetilde h_k\big),
\end{equation*}
and $\widetilde P_n$ is the orthoprojector onto  $\widetilde{\cal H}_n=\text{Span}\{h_{n+k-1}\}$, $k\ge1$.

Let $\widetilde x\in\widetilde{\cal H}_n$ be the maximizer of the Rayleigh quotient (\ref{Rayleigh}) over $\widetilde{\cal H}_n$. Just as in the first part, we obtain (\ref{J'=0}) for $k>n$.

Since $J(h_n)\ge \lambda_n$, we have $J(\widetilde x)\ge\lambda_n$. So, for any $k>n$ 
\begin{equation*}%\label{ak1}
\widetilde a_k:=|(\widetilde x,h_k)_{\cal H}|\le\frac {\lambda_n}{\lambda_n-\lambda_k}\,\sum\limits_{m=n}^{\infty} \widetilde a_m\cdot \frac c{(km)^{\frac {1+\delta}2}}.
\end{equation*}
%Without loss of generality we set $\widetilde a_n=1$. 
This implies 
%\begin{equation*}
%\widetilde A:=\sum\limits_{k=n}^{\infty} \frac {\widetilde a_k}{k^{\frac {1+\delta}2}}= \sum\limits_{k=n+1}^{\infty} \frac {\widetilde a_k}{k^{\frac {1+\delta}2}}+\frac 1{n^{\frac {1+\delta}2}}. 
%\end{equation*}
%Therefore, (\ref{ak1}) gives
\begin{equation}\label{A1}
\widetilde A:=\sum\limits_{k=n}^{\infty} \frac {\widetilde a_k}{k^{\frac {1+\delta}2}}\le \widetilde A\sum\limits_{k=n+1}^{\infty}\frac {\lambda_n}{\lambda_n-\lambda_k}\cdot \frac c{k^{1+\delta}}
+\frac 1{n^{\frac {1+\delta}2}}=:\widetilde A\widetilde {\mathfrak C}+\frac 1{n^{\frac {1+\delta}2}}.
\end{equation}

Notice that for $k>n$
\begin{equation*}
\frac {\lambda_n}{\lambda_n-\lambda_k}\le c\,\frac {(k/n)^B}{(k/n)^B-1}+\frac c{n^{\min\{1,B\}}}.
\end{equation*}
So, 
\begin{equation*}
\aligned
\widetilde {\mathfrak C}\le &\ \frac c{n^{1+\delta}}\sum\limits_{k=n+1}^{\infty} \frac {(k/n)^{B-1-\delta}}{(k/n)^B-1}+\frac c{n^{\min\{1,B\}}}\\
\le &\ \frac c{n^{\delta}}\int\limits_{1+\frac 1n}^{\infty} \frac {t^{B-1-\delta}}{t^B-1}\,dt+\frac c{n^{\min\{1,B\}}}\le \frac {c\log(n)}{n^{\min\{1,\delta,B\}}},
\endaligned
\end{equation*}
and for $n$ sufficiently large (\ref{A1}) gives $\widetilde A\le cn^{-\frac {1+\delta}2}$.

Now we calculate
\begin{equation*}
J(\widetilde x)\le\frac {\sum\limits_{k=n}^{\infty} \lambda_k \widetilde a_k^2}{\sum\limits_{k=n}^{\infty} \widetilde a_k^2-c\sum\limits_{k,m=n}^{\infty}
\dfrac {\widetilde a_k\widetilde a_m}{(km)^{\frac {1+\delta}2}}}\le \lambda_n(1+cn^{-(1+\delta)}),
\end{equation*}
and the statement again follows.
\hfill$\square$\medskip

\begin{Remark}
 The results of Lemma 1 and Theorem 1 hold true also in the case where the eigenvalues of the operator ${\cal K}$ are organized in two sequences
 \begin{equation*}
\lambda_n^{(1)}=\big((2n-1)a+b_1+O(n^{-\delta})\big)^{-B}, \qquad \lambda_n^{(2)}=\big(2na+b_2+O(n^{-\delta})\big)^{-B},
\end{equation*}
as $n\to\infty$. Such asymptotics is preserved under the same assumptions on the perturbation operator ${\cal B}$.
\end{Remark}

\end{document}